\documentclass[12pt]{article}

\usepackage{amssymb,amsfonts,amsmath,txfonts}
\usepackage{latexsym}

\newtheorem{theorem}{Theorem}

\newtheorem{proposition}[theorem]{Proposition}
\newtheorem{corollary}[theorem]{Corollary}

\newenvironment{proof}{\paragraph{\it Proof.}}{$\square$\vskip0.4cm}
\newenvironment{remark}{\paragraph{\it Remark.}}{\vskip0.4cm}

\newcommand{\nc}{\newcommand}
\nc{\beq}{\begin{eqnarray}}
\nc{\eeq}{\end{eqnarray}}
\nc{\bes}{\begin{eqnarray*}}
\nc{\ees}{\end{eqnarray*}}
\nc{\oper}[1]{\mathop{\mathchoice{\mbox{\rm #1}}{\mbox{\rm #1}}
{\mbox{\rm \scriptsize #1}}{\mbox{\rm \tiny #1}}}\nolimits}
\nc{\operlimits}[1]{\mathop{\mathchoice{\mbox{\rm #1}}{\mbox{\rm #1}}
{\mbox{\rm \scriptsize #1}}{\mbox{\rm \tiny #1}}}}

\nc{\gr}{{\mathfrak r}}
\nc{\gl}{{\mathfrak g}\mathfrak{l}}
\nc{\GL}{{\rm GL}}
\nc{\g}{{\mathfrak g}}
\renewcommand{\t}{{\mathfrak t}}
\nc{\T}{{\mathbb T}}
\nc{\F}{{\mathbb F}}
\newcommand{\vv}{{\mathbf v}}
\nc{\w}{{\mathbf w}}
\nc{\mF}{{\mathcal F}}
\nc{\mC}{{\mathcal C}}
\nc{\calM}{{\cal M}}
\nc{\calP}{{\cal P}}
\nc{\calT}{{\cal T}}
\nc{\N}{{\mathbb N}}
\nc{\Q}{{\cal Q}}
\nc{\J}{{\cal J}}
\nc{\higgsinfty}{${{\rm Higgs}_\infty}\ $}
\nc{\Qu}{{\mathbb Q}}
\nc{\frakC}{{\mathfrak C}}
\nc{\calC}{{\cal C}}
\nc{\calD}{{\cal D}}
\nc{\calV}{{\cal V}}
\nc{\calE}{{\cal E}}
\nc{\calB}{{\cal B}}
\nc{\C}{{\mathbb C}}
\nc{\R}{{\mathbb R}}
\nc{\calK}{{\cal K}}
\nc{\calR}{{\cal R}}
\nc{\Proj}{{\mathbb P}}
\nc{\Z}{{\mathbb Z}}
\nc{\Nat}{{\mathbb N}}
\nc{\E}{{\mathbb E}}
\nc{\U}{{\mathbb U}}
\nc{\dbar}{{\overline{\partial}}}
\nc{\K}{\mathbb K}
\nc{\Dirac}{{D\!\!\!\! \slash}} 
\nc{\Ha}{{\cal H}}
\nc{\Hy}{{\mathbb H}}
\nc{\V}{{\mathbb V}}
\nc{\M}{{\mathbb M}}
\nc{\I}{{\mathbb I}}
\renewcommand{\J}{{\mathbb J}}
\nc{\A}{{\mathcal A}}
\nc{\B}{{\mathcal B}}
\nc{\rmG}{{\rm G}}
\nc{\uPhi}{{\mathbf \Phi}}
\nc{\cM}{{{\overline\M}}}
\nc{\bM}{{\mathbf M}}
\nc{\calO}{{\cal O}}
\nc{\0}{{\mathbf 0}}
\nc{\bL}{{\mathbb L}}
\nc{\bP}{{\mathbb P}}
\nc{\bV}{{\mathbb V}}
\nc{\calL}{{\cal L}}
\nc{\cD}{{\cal D}}
\nc{\bd}{{\bar{d}}}
\nc{\bg}{{\bar{g}}}
\nc{\G}{{\cal G}}
\nc{\bG}{{\overline{\G}}}
\nc{\Pic}{\mathop{\rm Pic}\nolimits}
\nc{\coker}{\mathop{\rm coker}\nolimits}
\nc{\im}{\mathop{\rm im}\nolimits}
\nc{\rank}{\mathop{\rm rank}\nolimits}
\nc{\ch}{\mathop{\rm ch}\nolimits}
\nc{\td}{\mathop{\rm td}\nolimits}
\nc{\tr}{\mathop{\rm tr}\nolimits}
\nc{\pr}{\mathop{\rm pr}\nolimits}
\nc{\ad}{\mathop{\rm ad}\nolimits}
\nc{\Hom}{\mathop{\rm Hom}\nolimits}
\nc{\End}{\mathop{\rm End}\nolimits}
\nc{\Div}{\mathop{\rm Div}\nolimits}
\nc{\trace}{\mathop{\rm tr}\nolimits}
\nc{\sing}{\mathop{\rm Sing}\nolimits}
\nc{\const}{\oper{const.}}
\nc{\Coeff}{\operlimits{Coeff}}
\nc{\Res}{\operlimits{Res}}
\nc{\cE}{{\cal E}}
\nc{\compE}{{E\stackrel{\Phi}{\rightarrow}{E\otimes K}}}
\nc{\ucompE}{{\E_{\tM}\stackrel{\uPhi}{\rightarrow}
{\E_{\tM}\otimes\K}}}
\nc{\tD}{{\widetilde{D}}}
\nc{\tF}{{\widetilde{F}}}
\nc{\tM}{{\widetilde{\M}}}
\nc{\tN}{{\widetilde{\N}}}
\nc{\D}{{\mathbf D}}
\nc{\cL}{{\cal L}}
\nc{\eqed}{e^\circ_d}
\nc{\teqed}{\widetilde{e}^\circ_d}
\nc{\al}{\alpha}
\nc{\be}{\beta}
\nc{\ze}{\zeta}
\nc{\ga}{\gamma}
\nc{\la}{\lambda}
\nc{\La}{\Lambda}
\nc{\Ga}{\Gamma}
\nc{\si}{\sigma}
\nc{\xisi}{\xi^\Sigma}
\nc{\sisi}{\si^\Sigma}
\nc{\eqal}{\alpha^\circ}
\nc{\eqbe}{\beta^\circ}
\nc{\eqga}{\gamma^\circ}
\nc{\eqpsi}{\psi^\circ}
\nc{\eqze}{\zeta^\circ}
\nc{\eqeta}{\eta^\circ}

\nc{\bino}[2]{\mbox{\Large $#1 \choose #2$}}

\nc{\stack}[2]{
{\begin{array}{c}
\scriptstyle #1 \\ \scriptstyle #2 \end{array}} }

\begin{document}
\title{Betti numbers of holomorphic symplectic quotients via arithmetic Fourier transform}
\author{
 Tam\'as Hausel\footnote{Research supported by a Royal Society University Research Fellowship, an NSF grant DMS-0305505, an Alfred P. Sloan Research Fellowship and a Summer Research Assignement of the University of Texas at Austin.}
\\ {\it University of Oxford}
\\ {\it University of Texas at Austin}
\\{\tt hausel@maths.ox.ac.uk}}
\maketitle

\begin{abstract} A Fourier transform technique is introduced for counting the number of solutions of holomorphic moment map equations over a finite field. This in turn gives information
on Betti numbers of holomorphic symplectic quotients. 
As a consequence simple unified 
proofs are obtained  for formulas of Poincar\'e polynomials of toric
hyperk\"ahler varieties (recovering results of Bielawski-Dancer and
Hausel-Sturmfels), Poincar\'e polynomials of Hilbert schemes of points and twisted ADHM
spaces of instantons on $\C^2$ (recovering results of Nakajima-Yoshioka) and Poincar\'e polynomials of all Nakajima quiver varieties.
As an application, 
a proof of a conjecture of Kac on the number of absolutely
indecomposable representations of a quiver  is announced. 
\end{abstract}

Let $\K$ be a field, which will be either the complex numbers $\C$ or the finite field
$\F_q$ in this paper. Let 
$\rmG$ be a reductive algebraic group over $\K$, $\g$ its Lie algebra. Consider a representation $\rho: \rmG\to \GL(\V)$ of $\rmG$ 
on a $\K$-vector space $\V$, inducing the Lie algebra representation 
$\varrho: \g \to \gl(\V)$. 
This induces an action $\rho: \rmG\to \GL(\M)$ on $\M=\V\times \V^*$. 
 The vector
space $\M$ has a natural symplectic structure; defined by the natural pairing $\langle v,w\rangle =w(v)$, 
with $v\in \V$
and $w\in \V^*$. 
With respect to this symplectic form a moment map 
$$\mu: \V\times \V^* \to \g^*$$ 
of $\rho$ is given at $X\in \g$ by  \beq \label{momentmap} 
\langle \mu(v,w),X\rangle=\langle \varrho(X)v,w\rangle.
\eeq Let now $\xi\in (\g^*)^\rmG$ be 
a central element, then the holomorphic symplectic quotient is defined  by
the affine GIT quotient: $$\M/\!/\!/\!/\!_\xi \rmG:= (\mu^{-1}(\xi))/\!/\rmG,$$
which is the affine algebraic geometric version of the hyperk\"ahler
quotient construction of \cite{hitchin-etal}. In particular our varieties, additionally to the holomorphic symplectic structure, 
will carry a natural hyperk\"ahler metric, although the latter will 
not feature in what follows. 

Our main proposition counts rational points on the varieties $\mu^{-1}(\xi)$ over the finite fields $\F_q$, where
$q=p^r$ is a prime power. For convenience we will use the same letters $\V,\rmG,\g,\M,\xi$ for the corresponding vector
spaces, groups, Lie algebras and matrices over the finite field $\F_q$. We define the function 
$a_\varrho:\g\to {\mathbb N}\subset \C$ at $X\in \g$ as \beq \label{kernel} a_\varrho(X):=|\ker(\varrho(X))|, \eeq where we used the notation  $|S|$ for the 
number of elements in any set $S$.
In particular $a_\varrho(X)$ is always a power of $q$.  For an element $v\in V$ of any vector space
we define the characteristic function $\delta_v:V\to \C$ by $\delta_v(x)=0$ unless $x=v$ when $\delta_v(v)=1$.  We can now formulate
a generalization of the Fourier transform formula in \cite{hausel}:

\begin{proposition} The number of solutions of the equation  $\mu(v,w)=\xi$
 over the finite field $\F_q$ equals:  \bes \# \{(v,w)\in \M\,\, |\,\, \mu(v,w)= \xi\}&=&  |\g|^{-1/2} |\V| \mF(a_{\varrho
})(\xi)= |\g|^{-1} |\V| \sum_{X\in \g} a_{\varrho
}(X) \Psi(\langle X, \xi\rangle)\ees
\label{main}
\end{proposition} 

In order to explain the last two terms in the proposition above
we need to define
Fourier transforms \cite{letellier} of functions $f:\g\to \C$ on the finite Lie algebra $\g$, which here we think of as 
an abelian group with its additive structure. To define this fix $\Psi:\F_q\to \C^\times$ a non-trivial additive character,
and then we define the Fourier transform  $\mF(f):\g^*\to \C$ at a $Y\in \g^*$
$$\mF(f)(Y)= |\g|^{-1/2} \sum_{X\in \g}f(X) \Psi(\langle X, Y\rangle).$$

\begin{proof}  Using  two basic properties of Fourier transform: 
$$\mF(\mF(f))(X)=f(-X)$$  for $X\in \g$ and  \beq\label{basic2} \sum_{w\in V^*} \Psi(\langle v,w\rangle) = |V| \delta_0(v) \eeq  for $v\in V$ we get:
 \beq  \nonumber \# \{(v,w)\in \M\,\, |\,\, \mu(v,w)=\xi  \}&=& \nonumber \sum_{v\in \V} \sum_{w\in \V^*} 
\delta_\xi (\mu(v,w)) =  \sum_{v\in \V} \sum_{w\in \V^*} \mF(\mF(\delta_{\xi}))(-\mu(v,w))\\ &=&  \nonumber  \sum_{v\in \V} \sum_{w\in \V^*}
 \sum_{X\in \g}|\g|^{-1/2} \mF(\delta_{\xi})(X) \Psi(\langle X,-\mu(v,w)\rangle))\\ \nonumber &=& \sum_{v\in \V} 
\sum_{X\in \g}|\g|^{-1/2} \mF(\delta_{\xi})(X) \sum_{w\in \V^*} \Psi(-\langle\varrho(X)v,w\rangle )\\ 
 \nonumber &=& \sum_{v\in \V} 
\sum_{X\in \g}|\g|^{-1/2} \mF(\delta_{\xi})(X) |\V| \delta_0( \varrho(X)v)  \\ 
\nonumber &=&\sum_{X\in \g}|\g|^{-1/2} \mF(\delta_{\xi})(X) |\V| a_\varrho
(X)  \\ 
\nonumber &=&\sum_{X\in \g}|\g|^{-1} |\V| a_\varrho
(X) \sum_{Y\in \g^*} \delta_{\xi}(Y)  \Psi(\langle X,Y\rangle) \\ \nonumber &=& |\g|^{-1} |\V| \sum_{X\in \g} a_{\varrho
}(X) \Psi(\langle X, \xi\rangle) \\
\nonumber &=&  |\g|^{-1/2} |\V| \mF(a_{\varrho
})(\xi)\eeq 
\end{proof}

\newpage
\section{Affine toric hyperk\"ahler varieties}
We take $\rmG=\T^d\cong (\C^\times)^d$ a torus. A vector configuration
$A=(a_1,\dots,a_n):\Z^n\to \Z^d$ gives a representation 
$\rho_A:\T^d \to \T^n\subset \GL(\V)$, where $\V\cong \C^n$ is an $n$-dimensional vector space and $\T^n\subset \GL(\V)$ is a fixed maximal torus. The 
corresponding map on the Lie algebras is $\varrho_A: \t^d\to\t^n$. The 
holomorphic moment map of this action $\mu_A: \V\times \V^* \to (\t^d)^*$ is given by (\ref{momentmap})
which in this case takes the explicit form $$\mu_A (v,w) = \sum_{i=1}^n v_i w_i a_i.$$ We take a generic
$\xi\in (\t^d)^*$. The affine
toric hyperk\"ahler variety is then defined as the affine GIT quotient: 
$\calM(\xi,A)=\mu_A^{-1}(\xi)/\!/\T^d$.  In order to use our main result we need to 
determine $a_\varrho
(X)$. 
Note that the natural basis $e_1,\dots,e_n\in (\t^n)^*$ gives us a collection of 
hyperplanes $H_1,\dots,H_n$ in $\t^d$. Now for $X\in \t^d$
 we have that $a_\varrho
(X)=q^{ca(X)},$ where $ca(X)$ 
is the number of hyperplanes, which contain $X$. Finally we take the intersection lattice
$L(A)$ of this hyperplane arrangement; i.e. the set of all subspaces of
$\t^d$ which arise as the intersection of any collection of our hyperplanes; with partial
ordering given by containment. The generic choice of $\xi$ will
ensure that $\xi$ will not be trivial on any subspace in the lattice $L(A)$. Thus
for any subspace $V\in L(A)$, we have from (\ref{basic2}) that 
$\sum_{X\in V} \Psi(\langle X,\xi\rangle)=0$. 

Now we can use Proposition~\ref{main}.
If we perform the sum 
we get a combinatorial expression: $$\#(\calM(\xi,A))= \frac{q^{n-d}}{(q-1)^d} 
\sum_{X\in \t^d} a_\varrho
(X) \Psi(\langle X,\xi\rangle)=\frac{q^{n-d}}{(q-1)^d} 
\sum_{V\in L(A)} \mu_{L(A)}(V) q^{ca(V)}, $$ where $\mu_{L(A)}$ is the M\"obius function 
of the partially ordered set
$L(A)$, while $ca(V)$ is the number of coatoms, i.e. hyperplanes 
containing $V$. Because the count above is polynomial in $q$ 
and the mixed Hodge structure on $\calM(\xi,A)$ is pure we get that for the Poincar\'e polynomial we need to take the opposite 
of the count polynomial, i.e. substitute $q=1/t^2$ and multiply by $t^{4(n-d)}$. This
yields 

\begin{theorem} The Poincar\'e polynomial of the toric hyperk\"ahler variety is given
by $$P_t(\calM(\xi,A))=\frac{1}{t^{2d}(1-t^2)^d} \sum_{V\in L(A)} \mu_{L(A)}(V) (t^2)^{n-ca(V)}.$$
\end{theorem}
One can prove by a simple deletion contraction argument (and it also follows\footnote{I thank Ed Swartz for this reference.} 
from the second proof of Proposition 6.3.26 of \cite{white}) that for any matroid $\calM_\A$ and
its dual $\calM_\B$
$$\frac{1}{q^{d}(1-q)^d} \sum_{V\in L(A)} \mu_{L(A)}(V) (q)^{n-ca(V)}=h(\calM_\B),$$
where $$h(\calM_\B)=\sum_{i=0}^{n-d} h_i(\calM_\B) q^i$$ is the 
$h$-polynomial of the dual matroid $\calM_\B$. This way we recover 
a result of \cite{bielawski-dancer} and \cite{hausel-sturmfels}, for a more recent arithmetic proof see
\cite{proudfoot-webster} :
\begin{corollary}\label{toric} The  Poincar\'e polynomial of the toric hyperk\"ahler variety is given
by $$P_t(\calM(\xi,A))=h(\calM_B)(t^2),$$
where $B$ is a Gale dual vector configuration of $A$.
\end{corollary}

\section{Hilbert scheme of $n$-points on $\C^2$ and ADHM spaces}

Here $\rmG=\GL(V)$, where $V$ is an $n$-dimensional $\K$ vector space. 
We need three types of basic 
representations of $\rmG$. The adjoint representation
$\rho_{ad}: \GL(V)\to \GL(\gl(V))$, the defining representation 
$\rho_{def}=Id:\rmG\to \GL(V)$ and the trivial representations 
$\rho^k_{triv}=1:\rmG\to \GL(\K^k)$. Fix $k$ and $n$. Define $\V=\gl(V)\times V\otimes \K^k$,
$\M=\V\times \V^*$
and  $\rho:\rmG\to \GL(\V)$  
 by $\rho=\rho_{ad}\times \rho_{def}\otimes \rho^k_{triv}$. Then we 
take the central element $\xi=Id_V\in \gl(V)$ and define the twisted ADHM space  as  
$$\calM(n,k)=\M/\!/\!/\!/\!_\xi \rmG=\mu^{-1}(\xi)/\!/\rmG,$$  where 
$$\mu(A,B,I,J)=[A,B]+IJ,$$ with $A,B\in \gl(V)$, $I\in \Hom(\K^k, V)$ and
$J\in \Hom(V,\K^k)$. 

The  space $\calM(n,k)$ is empty when $k=0$ (the trace of a commutator
is always zero), 
diffeomorphic with the Hilbert scheme 
of $n$-points on $\C^2$, when $k=1$, and is the twisted version of the ADHM space \cite{atiyah-etal}
of $U(k)$ Yang-Mills instantons of charge $n$ on $\R^4$ (c.f. \cite{nakajima-book}). 
By our main Proposition~\ref{main}
 the number of solutions over $\K=\F_q$ of the equation $$[A,B]+IJ=Id_V$$ 
is the Fourier
transform on $\g$ of the function 
$a_\varrho
 (X)=|\ker(\varrho(X))|$. First we determine $a_\varrho
(X)$
for $X\in \g=\gl(V)$. By the definition of $\varrho$ we have $$\ker(\varrho(X))= \ker
(\varrho_{ad}(X)) \times \ker(\varrho_{def}) \otimes \K^k,$$ and so if 
$a_{\varrho_{ad}}(X)=|\ker
(\varrho_{ad}(X))|$ and $a_{\varrho_{def}}=|\ker(\varrho_{def})|$ then we have 
$$a_\varrho
(X)=a_{\varrho
_{ad}}(X) a^k_{\varrho_{def}}(X) .$$ This and Proposition~\ref{main} gives us \beq \nonumber \#(\calM(n,k))&=&\frac{1}{|\rmG|} \# \{ (v,w)\in \M | \mu (v,w)=\xi \}\\ \nonumber &=& \frac{|\V|}{|\g| |\rmG|} \sum_{X\in \g} a_{\varrho
}(X) \Psi(\langle X, \xi\rangle)\\ \nonumber &=&  \frac{|\V|}{|\g| |\rmG|} \sum_{X\in \g}a_{\varrho
_{ad}}(X) a^k_{\varrho_{def}}(X)  \Psi(\langle X, \xi\rangle).\eeq We will perform the 
sum adjoint orbit by adjoint orbit. The adjoint orbits of $\gl(n)$, according to their Jordan normal forms, fall into types,
labeled by $\calT(n)$, which stands for the set of all possible Jordan
normal forms of elements in $\gl(n)$. We denote by $\calT_{\!\! reg}(t)$ the types of the regular (i.e. non-singular) adjoint orbits, 
while $\calT_{\!\! nil}(s)=\calP(s)$ denotes the types of the nilpotent adjoint orbits, which
are just given by partitions of $s$. 
First we do the $k=0$ case where we know a priori, that the count should be $0$, because
the commutator of any two matrix is always trace-free thus cannot equal $\xi$ (for almost all $q$). 
Additionally, if we separate the nilpotent and regular parts of our adjoint orbits we get \beq \nonumber 0= \frac{1}{ |\rmG|} \sum_{X\in \g}a_{\varrho
_{ad}}(X)   
\Psi(\langle X, \xi\rangle) &=&   \sum_{n=s+t} 
\sum_{\lambda \in \calT_{\!\! nil}(s)} \frac{ |\frakC_\lambda|}{|C_\lambda|} 
\sum_{ \tau \in \calT_{\!\! reg}(t)} \frac{|\frakC_\tau|}{|C_\tau|}\Psi(\langle X_\tau, \xi\rangle), \eeq
where $C_\tau$ and respectively 
${\mathfrak C}_\tau$ denotes the centralizer 
of an element $X_\tau$ of $\g$ of type $\tau$ in the adjoint representation of $\rmG$, respectively 
$\g$ on $\g$. 

So if we define the generating serieses:
$$\Phi^0_{nil}(T)=1+\sum^\infty_{s=1} 
\sum_{\lambda \in \calT_{\!\! nil}(s)} \frac{|\frakC_\lambda|}{|C_\lambda|} T^s, $$ and
$$\Phi_{reg}(T)=1+\sum^\infty_{t=1}\sum_{\tau \in \calT_{\!\! reg}(t)} \frac{|\frakC_\tau|}{|C_\tau|}\Psi(\langle X_\tau, \xi\rangle)T^t,$$ then we have 
$$\Phi^0_{nil}(T) {\Phi_{reg}(T)}=1.$$ However $\Phi^0_{nil}$ is easy to calculate \cite{feit-fine}\footnote{I thank Fernando Rodriguez-Villegas for this reference.}:
\beq \label{phi0n} \Phi^0_{nil}(T)=\prod_{i=1}^{\infty} \prod_{j=1}^\infty \frac{1}{(1-T^iq^{1-j})},\eeq
thus we get \beq \label{phir} {\Phi_{reg}(T)}=\prod_{i=1}^{\infty} \prod_{j=1}^\infty {(1-T^iq^{1-j})}.\eeq
Now the general case is easy to deal with:
\beq \nonumber \frac{\#(\calM(n,k))}{q^{nk}}= \frac{1}{ |\rmG|} \sum_{X\in \g}a_{\varrho
_{ad}}(X)   
\Psi(\langle X, \xi\rangle) &=&   \sum_{n=s+t} 
\sum_{\lambda \in \calT_{\!\! nil}(s)} \frac{ |\frakC_\lambda|a^k_{\varrho
_{def}}(X_\lambda)}{|C_\lambda|} 
\sum_{ \tau \in \calT_{\!\! reg}(t)} \frac{|\frakC_\tau|}{|C_\tau|}\Psi(\langle X, \xi\rangle).
 \eeq Thus if we define 
the grand generating function by \beq \label{grand} \Phi^k(T)=1+\sum_{n=1}^\infty  
\#(\calM(n,k))\frac{T^n}{q^{kn}}\eeq and $$\Phi^k_{nil}(T)=1+\sum^\infty_{s=1} 
\sum_{\lambda \in \calT_{\!\! nil}(s)} \frac{|\frakC_\lambda| |\ker(X_\lambda)|^k}{|C_\lambda|} T^s, $$ then for the latter we get similarly to the argument for (\ref{phi0n}) in \cite{feit-fine} that $$\Phi^k_{nil}=  \Phi^k_{nil}(T)=\prod_{i=1}^{\infty} \prod_{j=1}^\infty \frac{1}{(1-T^iq^{k+1-j})}.$$ For the grand generating function then we get 
$$\Phi^k(T)=\Phi_{nil}^k(T)\Phi_{reg}(T)=\prod_{i=1}^{\infty} \prod_{l=1}^k \frac{1}{(1-T^iq^{l})}.$$ Because the mixed Hodge structure is pure, and this count is polynomial, this also
gives the compactly supported Poincar\'e polynomial. In order to get the ordinary 
Poincar\'e polynomial, we need to replace $q=1/t^2$ and multiply the $n$th  term in
(\ref{grand}) by $t^{4kn}$. This way we get

\begin{theorem} \label{adhm} The generating function of the Poincar\'e 
polynomials of the twisted ADHM spaces, are given by:
$$\sum_{n=0}^\infty P_t(\calM(k,n)) T^n=\prod_{i=1}^\infty \prod_{b=1}^k \frac{1}{(1-t^{2\left( k(i-1)+b-1\right)}T^i)}.$$\end{theorem} This result appeared as\footnote{I thank Bal\'azs Szendr\H oi for this reference.}  Corollary 3.10 in
\cite{nakajima-yoshioka}.

\section{Quiver varieties of Nakajima}

Here we recall the definition of the affine version of Nakajima's 
quiver varieties \cite{nakajima}. Let $Q=(\calV,\calE)$ be a quiver, i.e. an oriented 
graph on a finite set $\calV=\{1,\dots,n\}$ with $\calE\subset \calV\times \calV$ a finite set of oriented
(perhaps multiple and loop) edges. To each vertex $i$ of the graph we associate two finite dimensional $\K$ vector
spaces $V_i$ and $W_i$. We call $(\vv_1,\dots,\vv_n,\w_1,\dots, \w_n)=(\vv,\w)$ the dimension vector, where $\vv_i=\dim(V_i)$ and $\w_i=\dim(W_i)$. To this data we associate the grand vector space:
$$\V_{\vv,\w}=\bigoplus_{(i,j)\in \calE} \Hom(V_i,V_j) \oplus \bigoplus_{i\in \calV} \Hom (V_i,W_i),$$
the group and its Lie algebra $$\rmG_{\vv}=\varprod_{i\in \calV} \GL(V_i)$$ $$ \g_{\vv}=\bigoplus_{i\in \calV} \gl(V_i),$$ and the natural representation 
$$\rho_{\vv,\w}:\rmG_{\vv}\to \GL(\V_{\vv,\w}),$$ with derivative $$\varrho_{\vv,\w}:\g_{\vv}\to \gl(\V_{\vv,\w}).$$ The action is  from both left and right
on the first term, and from the left on the second.  

We now have $\rmG_{\vv}$ acting on $\M_{\vv,\w}=\V_{\vv,\w}\times \V_{\vv,\w}^*$ preserving
the symplectic form with moment map $\mu_{\vv,\w}:\V_{\vv,\w}\times\V_{\vv,\w}^*\to \g_{\vv}^*$ given by (\ref{momentmap}).
We take now $\xi_\vv=(Id_{V_1},\dots,Id_{V_n})\in (\g_{\vv}^*)^{\rmG_{\vv}}$, and define the affine Nakajima quiver variety \cite{nakajima} as 
$$\calM(\vv,\w)=\mu_{\vv,\w}^{-1}(\xi_\vv)/\!/\rmG_{\vv}.$$
Here we determine the Betti numbers of $\calM(\vv,\w)$ using our main Proposition~\ref{main}, by calculating the Fourier transform
of the function $a_{\varrho
_{\vv,\w}}$ given in (\ref{kernel}).

First we introduce,  for a dimension vector $\w\in \calV^\N$, the 
generating function 
$$\Phi_{nil}(\w)=
\sum_{\vv=(\vv_1,\dots,\vv_n)\in \calV^\N} \prod_{i\in V} T_i^{\vv_i} \sum_{\lambda^1\in 
\calT_{\!\! nil}(\vv_1)}\dots\sum_{{\lambda^n}\in \calT_{\!\! nil}(\vv_n)}  \frac{a_{\varrho
_{\vv,\w}}(X_{\lambda^1},\dots,X_{\lambda^n})}{|C_{\lambda^1}| \cdots |C_{\lambda^n}|}, 
$$ where $\calT_{\!\! nil}(s)$ is the set of types of nilpotent $s\times s$ 
matrices; where a type is given by a partition $\lambda\in \calP(s)$
of $s$, $X_\lambda$ denotes the typical $s\times s$ nilpotent matrix in $\gl(s)$ in Jordan form
of type $\lambda$, $C_\lambda$ is the centralizer of $X_\lambda$ 
under the adjoint action of $\GL(s)$ on $\gl(s)$.
We also introduce the generating function 
$$\Phi_{reg} = \sum_{\vv=({\vv}_1,\dots,\vv_n)\in \calV^\N} \prod_{i\in V} T_i^{\vv_i} \sum_{\tau_1\in 
\calT_{\!\! reg}(\vv_1)}\dots\sum_{{\tau_n}\in \calT_{\!\! reg}(\vv_n)}  \frac{a_{\varrho
_{\vv,\0}}(X_{\tau_1},\dots,X_{\tau_n})}{|C_{\tau_1}| \cdots |C_{\tau_n}|}\Psi(\langle X_{\tau}, \xi_\vv\rangle), 
$$ where $\calT_{\!\! reg}(t)$ is the set of types $\tau$, i.e. Jordan normal forms, of a regular $t\times t$ matrix $X_\tau$
in $\gl(t)$, $C_\tau\subset \GL(t)$ its centralizer under the adjoint action. Note also that for a regular element $X\in \g_{\vv}$, $a_{\varrho
_{\vv,\w}}(X)=a_{\varrho
_{\vv,\0}}(X)$ does not depend on $\w\in \calV^\N$. 

Now  we introduce for $\w\in \calV^\N$ 
the grand generating function \beq \label{grandquiver}\Phi(\w)=
\sum_{\vv\in \calV^\N} \#(\calM(\vv,\w))  \frac{|\g_{\vv}|}{|\V_{\vv,\w}|}T^\vv.\eeq
As in the previous section, our main Proposition~\ref{main} implies \beq \label{functional} 
\Phi(\w)=\Phi_{nil}(\w) \Phi_{reg}.\eeq Finally we note, that when $\w=\0$ we have $\varrho_{\vv,\0}(\xi_\vv^*)=0$, where $\xi_\vv^*=(Id_{V_1},\dots,Id_{V_n})\in \g$, thus
by (\ref{momentmap})
$\langle \mu_{\vv,\0}(v,\0), \xi_\vv^*\rangle=0$. Because $\langle \xi_\vv,\xi_\vv^*\rangle=\sum \vv_i$, the equation $\mu_{\vv,\0}(v,w)=\xi_\vv$ 
has no solutions (for almost all $q$). This way we get that $\Phi(\0)=1$ and so (\ref{functional}) yields  $\Phi_{reg}=\frac{1}{\Phi_{nil}(\0)}$, giving the result $$\Phi(\w)=\frac{\Phi_{nil}(\w)}{\Phi_{nil}(\0)}.$$ Therefore it is enough to understand
$\Phi_{nil}(\w),$ which reduces to a simple linear algebra problem
of determining $a_{\varrho
_{\vv,\w}}(X_{\lambda^1},\dots,X_{\lambda^n}).$ Putting together everything yields the following:
\begin{theorem} Let $Q=(\calV,\calE)$ be a quiver, with $\calV=\{1,\dots,n\}$ and $\calE\subset \calV\times \calV$, with possibly multiple edges and loops. Fix a dimension vector  $\w\in \N^\calV$. The Poincar\'e polynomials $P_t(\calM(\vv,\w))$ of the corresponding Nakajima quiver varieties are given by the generating function:
\beq\label{quivergenerate}  \sum_{\vv\in \N^\calV} P_t(\calM(\vv,\w)) t^{-d(\vv,\w)}T^\vv=
\frac{\sum_{\vv\in \N^\calV} T^\vv \sum_{\lambda^1\in 
\calP(\vv_1)}\dots \sum_{\lambda^n\in \calP(\vv_n) } \frac{
\left( \prod_{(i,j)\in \calE} t^{-2n(\lambda^i,\lambda^j)}\right)\left(\prod_{i\in \calV} t^{-2n(\lambda^i,(1^{\w_i}))}\right) }{\prod_{i\in \calV} \left(t^{-2n(\lambda^i,\lambda^i))}\prod_k \prod_{j=1}^{m_k(\lambda^i)} (1-t^{2j}) \right)}}{\sum_{\vv\in \N^\calV} T^\vv \sum_{\lambda^1\in 
\calP(\vv_1)}\dots\sum_{{\lambda^n}\in \calP(\vv_n)}  \frac{\prod_{(i,j)\in \calE} t^{-2n(\lambda^i,\lambda^j))}}{\prod_{i\in \calV} \left(t^{-2n(\lambda^i,\lambda^i))}\prod_k \prod_{j=1}^{m_k(\lambda^i)} (1-t^{2j})\right) }},\eeq where $d(\vv,\w)=2{{\sum_{(i,j)\in \calE} \vv_i\vv_j}+2{\sum_{i\in \calV} \vv_i(\w_i-\vv_i)}}$ is the dimension of $\calM(\vv,\w)$ and $T^\vv=\prod_{i\in \calV} T_i^{\vv_i}$.
 $\calP(s)$ stands for the set of partitions\footnote{ The notation
 for partitions is
that of \cite{macdonald}.}
 of $s\in \N$.  For two partitions $\lambda=(\lambda_1,\dots,\lambda_l)\in \calP(s)$ and $\mu=(\mu_1,\dots,\mu_m)\in \calP(s)$ we define 
$n(\lambda,\mu)=\sum_{i,j} \min(\lambda_i,\mu_j)$, and if we write $\lambda=(1^{m_1(\lambda)},2^{m_2(\lambda)},\dots)\in \calP(s)$, then we can define $l(\lambda)=\sum m_i(\lambda)=l$ the number of parts in $\lambda$ .  With this notation $n(\lambda^i,(1^{\w_i}))=\w_i l(\lambda^i)$ in the above formula.
\end{theorem}
\begin{remark} This single formula encompasses  a surprising amount of combinatorics and representation theory.  When $\vv=(1,\dots,1)$ the Nakajima quiver variety
is a toric hyperk\"ahler variety, thus (\ref{quivergenerate}) 
gives a new formula for its Poincar\'e polynomial, which was given in Corollary~\ref{toric}.
If additionally $\w=(1,0,\dots,0)$ then $\calM(\vv,\w)$ is the 
 toric quiver variety of \cite{hausel-sturmfels}. Therefore its Poincar\'e polynomial, which is the reliability polynomial\footnote{Incidentally, the reliability polynomial measures the probability of the graph 
 remaining connected when each edge has the same probability of failure; a concept heavily used in the study of reliability of computer networks \cite{colbourn}.}  of the graph underlying the quiver \cite{hausel-sturmfels} can also be read off from the above formula (\ref{quivergenerate}).

When the quiver is just a single loop on one vertex, our formula (\ref{quivergenerate}) reproduces Theorem~\ref{adhm}. When the quiver is of type $A_n$ Nakajima \cite{nakajima-anquiver} showed, that the Poincar\'e polynomials of the quiver variety are related to the combinatorics
of Young-tableaux, while in the general $ADE$ case, Lusztig \cite{lusztig}
conjectured a formula for the Poincar\'e polynomial, in terms of 
formulae arising in the representation theory of quantum groups. 
When the quiver is star-shaped recent work in \cite{hausel,hausel-etal} calculates these Poincar\'e polynomials using
the character theory of reductive Lie algebras over finite fields 
\cite{letellier}, 
and arrives at formulas determined by the Hall-Littlewood symmetric
functions \cite{macdonald}, which arose in the context of \cite{hausel,hausel-etal} as the pure part of Macdonald
symmetric polynomials \cite{macdonald}. In the case when the quiver has no loops, Nakajima \cite{nakajima-tq} gives a combinatorial algorithm for all Betti numbers of quiver varieties, motivated by
the representation theory of quantum loop algebras.  Finally, through the paper \cite{crawley-boevey-etal} of Crawley-Boevey and Van den Bergh, Poincar\'e polynomials
of quiver varieties are related 
to the number
of absolutely indecomposable representations of quivers in the work of
Kac \cite{kac} ; which were eventually 
completely determined by Hua \cite{hua}. 

In particular, formula (\ref{quivergenerate}), when combined with results in \cite{hua}, \cite{nakajima} and the Weyl-Kac character formula in the representation theory of Kac-Moody algebras \cite{kac0} yields\footnote{I thank Hiraku Nakajima
for suggesting this possibility.} a simple proof of Conjecture 1 of Kac 
\cite{kac}. Consequently, formula (\ref{quivergenerate}) can be viewed
as a $q$-deformation of the  Weyl-Kac character formula \cite{kac0}.
  
A detailed  study of the above generating function (\ref{quivergenerate}), its relationship to the wide variety of examples
mentioned above and details of the proofs of the results of this paper  will appear 
elsewhere. 

\end{remark}

\end{document}